\documentclass[10pt]{amsart}
\usepackage{epsfig, psfrag}
\usepackage[dvips]{xypic}  
\usepackage{bbm}

\xyoption{curve}                

\newcommand{\Z}{\mathbb{Z}}
\newcommand{\R}{\mathbb{R}}
\newcommand{\wt}{\widetilde}
\newcommand{\la}{\langle}
\newcommand{\ra}{\rangle}

\newcommand{\lie}{{{\mathcal L}{\it ie\/}}}
\newcommand{\ass}{{{\mathcal A}{\it ss\/}}}
\newcommand{\ndo}{{{\mathcal E}{\it\/nd\/}}}
\newcommand{\dis}{{{\mathcal D}{\it isks\/}}}
\newcommand{\com}{{{\mathcal C}{\it om\/}}}

\newcommand{\iopd}{{{\mathcal S}{\it iop\/}^{d}}}
\newcommand{\poidn}{{{\mathcal P}{\it ois\/}^{d}(n)}}

\newcommand{\poid}{{{\mathcal P}{\it ois\/}^{d}}}
\newcommand{\n}{\mathbf{n}}
\newcommand{\ka}{\mathbbm{k}}
\newcommand{\Op}{{\mathcal O}}  
\newcommand{\mC}{{\mathcal{C}}}
\newcommand{\one}{{\mathbbm{1}}}
\newcommand{\Co}{{\rm Conf}}
\newcommand{\col}{{\rm col}}
\newcommand{\dlog}{{\rm dlog}\;}

\newtheorem{theorem}{Theorem}[section]
\newtheorem{corollary}[theorem]{Corollary}
\newtheorem{lemma}[theorem]{Lemma}
\newtheorem{proposition}[theorem]{Proposition}

\theoremstyle{definition}                          
\newtheorem{definition}[theorem]{Definition}  
\newtheorem{examples}[theorem]{Examples}
\theoremstyle{remark}

\newcommand{\refT}[1]{Theorem~\ref{T:#1}}
\newcommand{\refC}[1]{Corollary~\ref{C:#1}}
\newcommand{\refP}[1]{Proposition~\ref{P:#1}}
\newcommand{\refD}[1]{Definition~\ref{D:#1}}
\newcommand{\refL}[1]{Lemma~\ref{L:#1}}

\newcommand{\refF}[1]{Figure~\ref{F:#1}}
\newcommand{\refS}[1]{Section~\ref{S:#1}}

\newcommand{\tree}[2]{ \ensuremath{  
 \begin{xy}                          %
   (0,1); (1,0)**\dir{-};            
   (2,1)**\dir{-},                   
   (1,-1); (1,0)**\dir{-},           
   (0,2.2)*{\scriptstyle #1},        %
   (2,2.2)*{\scriptstyle #2},        %
 \end{xy}  } } 

\newcommand{\linep}[2]{ \ensuremath{  %
 \begin{xy}                           
  (-1,-1)*+UR{\scriptstyle #1}="a",    
  (4,2)*+UR{\scriptstyle #2}="b",     
  "a";"b"**\dir{-}?>*\dir{>},         
 \end{xy}                             %
} }

\begin{document}

\title{The (non-equivariant) homology of the little disks operad}
\author{Dev P. Sinha}
\address{Department of Mathematics, University of Oregon, Eugene OR 97403, USA}
\maketitle

\section{Introduction}

This expository paper aims to be a gentle introduction to the topology of configuration
spaces, or equivalently spaces of little disks.  The pantheon of topological spaces
which beginning graduate students see is  limited -- spheres, projective spaces,
products of such, perhaps some spaces such as Lie groups, 
Grassmannians or knot complements.  We would like
for Euclidean configuration spaces to be added to this list.
We aim for this article 
to be appropriate  for someone who knows only basic homology and cohomology theory.  The
one exception to this rule will be the light use of a spectral sequence argument,
for an upper bound. 

Any space from the pantheon  has rich associated combinatorial and algebraic structure.
For example, the relationship between cohomology of
projective spaces and Grassmannians is encoded by the structure of symmetric polynomials.
In the case of configuration spaces, we are led to study graphs, trees, Jacobi and Arnold
identities, and ultimately the Poisson operad.   One goal of this paper
is to explain the topology which leads to the configuration pairing between
graphs and trees, developed 
purely combinatorially in \cite{Sinh06.2}.  That this pairing arises as that between canonical
spanning sets for  homology and cohomology
of configuration spaces is a new result.  Another goal is to prepare a reader
for further study of the theory of operads by giving a thorough understanding of the
disks operad from topology, the Poisson operad from algebra, and the fact that they are related
through homology.  We also bring in recently developed ingredients such as canonical 
compactifications of configuration spaces and submanifolds defined by 
collinearities.  These  new results and points of view,
and our elementary development, differentiate this paper from expositions 
such as \cite{CohF95, CLM76, FaHu01, CohF73, Arno69}.

The plan of the paper is as follows.   First in \refS{homology} 
we associate  a  class in the homology of configuration spaces 
to any forest, as the fundamental class of a submanifold homeomorphic to a torus,
and then develop relations between such classes.   Then in \refS{cohomology} we associate a cohomology class
 which is pulled back from a map to a torus to any graph.  \refS{pairing} gives the main new results, identifying the
evaluation of our graphical cohomology classes on the forest homology classes with a combinatorially defined
pairing between graphs and trees.  This pairing is useful in a number of contexts, for example in simultaneously 
understanding free Lie algebras and coalgebras \cite{SiWa06}, but is not widely known.   In \refS{operads} we 
give an informal ``examples-first'' development of operads, complementary to others in this volume, in 
order to be self-contained.   In \refS{diskshomology}, in particular \refT{main2}, 
we prove the well-known result that the homology of the little
disks operad is the graded Poisson operad.  Instead of the usual practice of
waiving our hands at the operad structure maps, 
we are able to provide a complete argument by arguing on cohomology instead.
At the end of most sections we give some (incomplete) historical notes.

The author would like to thank Ben Walter
for useful discussions, thank his students for feedback, and thank the organizers
of both the Banff graduate workshop on homotopy theory in 2005 and
the Luminy graduate workshop on operads in 2009 for their
encouragement in making this material accessible.

\tableofcontents

\section{Homology generators of configuration spaces}\label{S:homology}

In this section we construct homology classes for configuration
spaces,  all represented by submanifolds homeomorphic
to tori.  
The term ``configuration space'' is used in different ways by different subfields.
We use the term as is standard in algebraic topology, as the space of distinct
labeled points in some ambient space. 

\begin{definition}
The configuration space of $n$ distinct points in a space $X$, denoted $\Co_n(X)$,
is the subspace of the product $X^{\times n}$ defined as follows 
$$\{(x_1, \ldots, x_n) \in X^{\times n} | \, x_i \neq x_j \, {\rm if} \, i \neq j \}.$$
\end{definition}

We will sometimes abbreviate $(x_1, \cdots, x_n)$ as ${\bf x}$.  
We focus on the case in which $X$ is a Euclidean space $\R^d$.  This
configuration space models all possible simultaneous positions of $n$ 
distinct planets or particles.   
This space as a whole may be visualized through linear algebra, starting
with the ambient Euclidean space $\R^{nd}$ and removing the hyperplanes
where some $x_i = x_j$.  Indeed, the Euclidean configuration spaces are 
special important cases of complements of hyperplane arrangements.

Our strategy to compute the homology and cohomology of these spaces
is to ``just get our hands on things.''  When $n=2$
we have the following.

\begin{proposition}\label{P:n=2}
The configuration space $\Co_2(\R^d)$ is homotopy equivalent to $S^{d-1}$.
Thus the homology of $\Co_2(\R^d)$ is free, rank one in dimensions $0$ and $d-1$,
and zero otherwise.
\end{proposition}

\begin{proof}
We include $S^{d-1}$ into $\Co_2(\R^d)$ as a deformation retract.  With an
eye towards generalization, define the subspace $P_{\tree{1}{2}}$ of $\Co_2(\R^d)$
as $\{(x_1, x_2) \, | \, x_1 = -x_2 \, {\rm and} \, |x_i| = 1\}.$  The deformation
retract onto this subspace sends $(x_1, x_2)$ to 
$\left(\frac{x_1 - m}{|x_1 - m|}, \frac{x_2 - m}{|x_2 - m|}\right)$, where $m = \frac{x_1 + x_2}{2}$.  
The homotopy between this retraction and the identity map
is given by a straight-line homotopy. 
\end{proof}

We also deduce that the generating cycle in $H_{d-1}\left(\Co_2(\R^d)\right)$ is the image of
the fundamental class of the sphere by the map which sends $v \in S^{d-1}$
to $(v, -v)$, parameterizing $P_{\tree{1}{2}}$.  
Dually, $H^{d-1}\left(\Co_2(\R^d)\right)$ is pulled back from the sphere
by the given retraction.  More geometrically, a generator of cohomology
is Lefshetz dual to the submanifold $(x_1, x_2)$ such that $\frac{x_1 - x_{2}}{|x_1 - x_{2}|}$
is say the north pole 
$S^{d-1}$. This cohomology thus evaluates on some $d-1$ dimensional cycle by
counting with signs the number of configurations parameterized by that cycle
for which ``$x_1$ lies over $x_2$.'' 


 \medskip

For the general case the language of solar systems is suggestive, as Fred Cohen likes to point out.
In the $n=2$ case, the fundamental cycle had the ``planets'' $x_1$
and $x_2$ ``orbiting'' their center of mass.  For $n>2$ we can build further
cycles by having that ``system''  orbit the common center of mass
with some other planet or system of planets.  Each time we build a  system, 
it is possible (if there are more planets around) to put
that in an orbit with another system to create a more complicated one.  
Such systems, which are more difficult to formalize
than to visualize (see Figure~1), are naturally
indexed by trees. 

\begin{definition}\label{D:trees}
\begin{enumerate}
\item An $S$-tree is an isotopy class of acyclic graph whose vertices
are either trivalent or univalent, with a 
distinguished univalent vertex called the root, embedded in the 
upper half-plane with the root at the origin  -- for example $T =
\begin{xy}
  (1.5,1.5); (3,3)**\dir{-}, 
  (0,3); (3,0)**\dir{-}; 
  (7.5,4.5)**\dir{-},        
  (3,0); (3,-1.5)**\dir{-}, 
  (3,6); (6,3)**\dir{-},
  (6,6); (4.5,4.5)**\dir{-},
  (0,4.2)*{\scriptstyle 2},
  (3,4.2)*{\scriptstyle 6},
  (3,7.2)*{\scriptstyle 1},
  (6,7.2)*{\scriptstyle 7},
  (7.5,5.7)*{\scriptstyle 3},
\end{xy}$.  
Univalent vertices other than the root are called leaves, and they
are labeled by a subset $S$ of some set $\n = \{ 1, \ldots, n \}$.
Trivalent vertices are also called internal vertices.
\item The height of a vertex in an $S$-tree, denoted $h(v)$, is the number of edges between 
$v$ and  the root.  Edges which connect a vertex to higher vertices are called outgoing.
\item To define a subtree of $T$, take some vertex $v$ and all of the 
vertices and edges above it. Restrict the ambient embedding in the upper half plane,
and add a root edge from $v$ to the origin, to obtain a tree we call $T_v$.
Moreover, let $T_v^L$ be the subtree associated to the left vertex over $v$,
and similarly $T_v^R$ be the right subtree over $v$.
\item We say that $v$ is above or over $w$ if $w$ lies in the shortest path from
$v$ to the root.  Define a total order on the vertices of $T$ so that $v < w$ if $v$ 
lies over the left outgoing edge of $w$ and $v > w$ if it lies over the right 
outgoing edge of $w$.  This total ordering can be realized as a left
to right ordering of an appropriate planar embedding.

\end{enumerate}
\end{definition}

We now define the ``centers of mass'' for our systems and sub-systems.

\begin{definition}
The center $c({\bf x}, T)$ of a configuration ${\bf x}$ with respect to a tree
$T$ is defined inductively by $c({\bf x}, T_{v}) = \frac{1}{2}\left( c({\bf x}, {T_v^L}) +
c({\bf x}, {T_v^R}) \right)$, if $T$ has at least one internal vertex.  If $T$ consists
of only a leaf labeled by $i$, then $c({\bf x}, T) = x_i$.
\end{definition}

Finally, we can define the systems as ones where planets in a (sub)system
are of a prescribed distance from the center of mass.  Fix (for the moment)
an $\varepsilon < \frac{1}{3}$.

\begin{definition}\label{D:PT}
Given an $S$-tree $T$, the (planetary system) $P_T$ is the submanifold
of all ${\bf x} = (x_1, \cdots, x_n)$ such that:
\begin{enumerate}
\item $c({\bf x}, T) = 0$. \label{p1}
\item For any vertex $v$ of $T$, 
$ d\left(c({\bf x}, {T_v^L}), c({\bf x}, {T_v})\right) = \varepsilon^{h(v)} = 
 d\left(c({\bf x}, {T_v}), c({\bf x}, {T_v^R})\right),$ where $d$ is the standard
 Euclidean distance function. \label{p2}
\item If $i \notin S$, $x_i$ is fixed as some point ``at infinity.'' \label{p3}
\end{enumerate}
\end{definition}

\begin{figure}[ht]\label{F:1}
\psfrag{J}{$x_2$}
\psfrag{K}{$x_6$}
\psfrag{H}{$x_3$}
\psfrag{E}{$x_1$}
\psfrag{F}{$x_7$}
\psfrag{R}{$x_4$}
\psfrag{U}{$x_5$}
$$\includegraphics[width=12cm]{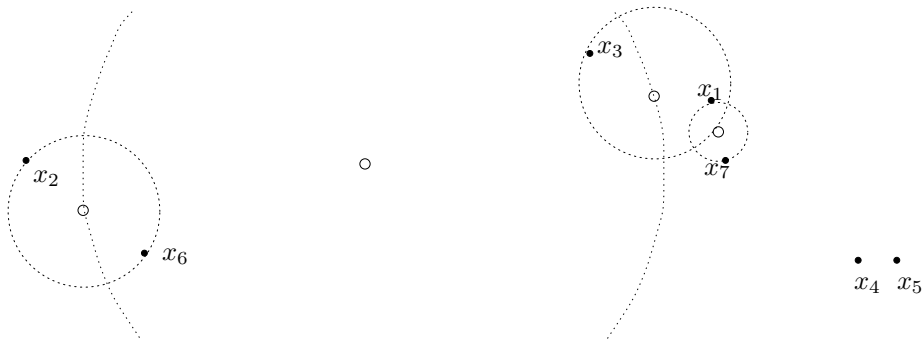}$$
\caption{An illustration of $P_T$}  \end{figure}

We picture these submanifolds as in Figure~1, which illustrates the case of 
$T =
\begin{xy}
  (1.5,1.5); (3,3)**\dir{-}, 
  (0,3); (3,0)**\dir{-}; 
  (7.5,4.5)**\dir{-},        
  (3,0); (3,-1.5)**\dir{-}, 
  (3,6); (6,3)**\dir{-},
  (6,6); (4.5,4.5)**\dir{-},
  (0,4.2)*{\scriptstyle 2},
  (3,4.2)*{\scriptstyle 6},
  (3,7.2)*{\scriptstyle 1},
  (6,7.2)*{\scriptstyle 7},
  (7.5,5.7)*{\scriptstyle 3},
\end{xy}$

One configuration in this submanifold is illustrated by the $\bullet$, which are labeled by
 $x_i$.  The rest of the family is indicated by drawing some of 
 the circular orbits where
 points in these configurations occur.  The centers of this configuration,
namely the points $c({\bf x}, {T_v})$, are indicated by $\circ$.
 In any configuration in $P_T$, the points $x_4$ and $x_5$ occur 
 where they are indicated.

We will use $P_{T}$ to define a homology class but to do so integrally, rather
than only with $\Z/2$ coefficients, it must be oriented.  We orient $P_{T}$
by parametrizing it through a map from a torus.

\begin{definition}\label{D:param}
By abuse of notation, let $P_{T} : (S^{d-1})^{\times|T|} \to \Co_{n}(\R^{d})$,
where $|T|$ is the number of internal vertices of $T$, send
$(u_{v_{1}}, \ldots, u_{v_{|T|}})$ to $(x_{1}, \ldots, x_{n})$, where
$$x_{i} = \sum_{v_{j} {\; \rm below \; leaf} \; i} \pm \varepsilon^{h_{j}} u_{v_{j}}.$$
Here the sum is taken over all vertices $v_{j}$ which lie on the path
from the leaf labeled by $i$ to the root vertex, and $h_{j}$ is the
height of $v_{j}$.  The sign $\pm$ is $+1$ if the path from 
the leaf $i$ to the root goes through the left edge of $v_{j}$ and
$-1$ if that path goes through the right edge of $v_{j}$.
\end{definition}

We may now orient $P_{T}$ by fixing an orientation of the sphere
and using the product orientation for $ (S^{d-1})^{\times|T|}$.
We will call the resulting homology class simply 
$T \in H_{|T|(d-1)}(\Co_{n}\left(\R^{d})\right)$.  Note that by its definition, $T$ is 
in the image of the map from the oriented bordism of  $\Co_{n}(\R^{d})$.
In fact, because spheres are stably framed it is in the 
image of the map from framed bordism, or equivalently stable homotopy.

The relations
between these homology classes represent a fundamental blending of geometry and algebra.

\begin{proposition}\label{P:relations}
The classes in $H_{*}(\Co(\R^{d}))$ given by trees satisfy the following relations:
\begin{align*}
\text{(anti-symmetry)} \qquad 
 & \qquad
\begin{xy}
   (0,1.5); (1.5,0)**\dir{-};          
   (3,1.5)**\dir{-},               
   (1.5,-1.5); (1.5,0)**\dir{-},       
   (-.4,2.7)*{\scriptstyle T_1},   
   (3.8,2.7)*{\scriptstyle T_2},   
   (1.5,-2.7)*{\scriptstyle R}
\end{xy}\  -
  (-1)^{d + |T_{1}||T_{2}|(d-1)} \begin{xy}
   (0,1.5); (1.5,0)**\dir{-};          
   (3,1.5)**\dir{-},               
   (1.5,-1.5); (1.5,0)**\dir{-},       
   (-.4,2.7)*{\scriptstyle T_2},   
   (3.8,2.7)*{\scriptstyle T_1},    
   (1.5,-2.7)*{\scriptstyle R}
\end{xy}  =\ 0\\
 \text{(Jacobi)} \qquad 
 & \qquad 
 \begin{xy}   
   (1.5,1.5); (3,3)**\dir{-}, 
   (0,3); (3,0)**\dir{-};
   (6,3)**\dir{-},   
   (3,0); (3,-1.5)**\dir{-}, 
   (-.4,4.2)*{\scriptstyle T_1}, 
   (3.2,4.2)*{\scriptstyle T_2},
   (6.8,4.2)*{\scriptstyle T_3}, 
   (3,-2.7)*{\scriptstyle R}
 \end{xy}  \ + \ \begin{xy}   
   (1.5,1.5); (3,3)**\dir{-}, 
   (0,3); (3,0)**\dir{-};
   (6,3)**\dir{-},   
   (3,0); (3,-1.5)**\dir{-}, 
   (-.4,4.2)*{\scriptstyle T_2}, 
   (3.2,4.2)*{\scriptstyle T_3},
   (6.8,4.2)*{\scriptstyle T_1}, 
   (3,-2.7)*{\scriptstyle R}
 \end{xy} \ + \ \begin{xy}   
   (1.5,1.5); (3,3)**\dir{-}, 
   (0,3); (3,0)**\dir{-};
   (6,3)**\dir{-},   
   (3,0); (3,-1.5)**\dir{-}, 
   (-.4,4.2)*{\scriptstyle T_3}, 
   (3.2,4.2)*{\scriptstyle T_1},
   (6.8,4.2)*{\scriptstyle T_2}, 
   (3,-2.7)*{\scriptstyle R}
 \end{xy} 
 \ =\  0
\end{align*}
where $R$, $T_1$, $T_2$, and $T_3$ stand for arbitrary (possibly trivial) 
subtrees which are not modified in these operations, and $|T_{i}|$ denotes
the number of internal vertices of $T_{i}$.
\end{proposition}

We call this relation the Jacobi identity because of the standard translation
between $S$-trees and bracket expressions, under which this becomes
$[[T_{1}, T_{2}], T_{3}] + [[T_{2}, T_{3}], T_{1}] + [[T_{3}, T_{1}], T_{2}] = 0$.

\begin{proof}
The anti-symmetry relation follows because the submanifolds defined
by these two trees are the same, and their parametrizations differ only by
the antipodal map on one factor of $S^{d-1}$ from \refD{param} and
reordering by moving the factors labeled by vertices
of $T_{1}$ after those of $T_{2}$

The Jacobi identity follows from the existence of Jacobi manifolds who
bound the submanifolds in that relation.  Letting $T$ be the first tree pictured in the Jacobi
identity above, consider the submanifold of $\Co_{n}(\R^{d})$
defined by conditions (\ref{p1}) and (\ref{p3}) from \refD{PT}, as well as condition~(\ref{p2}) 
for vertices  internal to $T_{1}$, $T_{2}$, $T_{3}$ or $R$.  For the remaining
two vertices we replace condition~(\ref{p2}) by
\begin{enumerate}
 \item[(2.1)] $\sum_{i, j \in \{1, 2, 3\}} 
 d\left(c({\bf x}, {T_{i}}), c({\bf x}, {T_{j}}) \right) = 4\varepsilon^{h} + 2\varepsilon^{h+1}$,
 where $h$ is the height of the internal vertex immediately above the
 subtree  $R$. \label{j3}
 \item[(2.2)] $d\left(c({\bf x}, {T_{i}}), c({\bf x}, {T_{j}}) \right)  \geq 2 \varepsilon^{h+1}$, 
 where $i,j \in \{1, 2, 3\}$. \label{j4}
\end{enumerate}
Condition~(2.1) fixes the perimeter of the triangle with vertices at the centers
of the sub-configurations associated to $T_{1}$, $T_{2}$, and $T_{3}$, and condition~(2.5)
says that triangle must have a minimum side length of at least $2\varepsilon^{h+1}$.

These conditions determine a submanifold $J$
with boundary, whose boundary is where one of the  three distances
$d\left(c({\bf x}, {T_{i}}), c({\bf x}, {T_{j}}) \right)$ is equal to $2 \varepsilon^{h+1}$.  There
are thus three boundary components.
By condition~(2.1), when some
$d\left(c({\bf x}, {T_{i}}), c({\bf x}, {T_{j}}) \right) = 2 \varepsilon^{h+1}$
the remaining center must be 
distance roughly $2 \varepsilon^{h}$ from the other two.  So these 
components of $\partial J$ are close to being the 
submanifolds $P_{T}$ for the $T$ which occur in the
Jacobi identity -- see Figure~\ref{F3}.
We get exactly those manifolds by replacing condition~(2.1) by
\begin{enumerate}
\item[(2.1')] $d\left(c({\bf x}, {T_{i}}), c({\bf x}, {T_{j}}) \right) = f({\bf x})$.
\end{enumerate}
Here $f({\bf x})$ is an interpolation function. Its value is 
$4\varepsilon^{h} + 2\varepsilon^{h+1}$ when the $c{{\bf x}, T_{i}}$ form a nearly equilateral
triangle.  When the configuration in question is in $P_{T}$, its value is 
the total length of the triangle with vertices
at $c({\bf x},{T_{i}})$, namely
$$ \sqrt{4 \varepsilon^{2h} + \varepsilon^{4h} - 4 \varepsilon^{3h} \cos(\theta)} + 
\sqrt{4 \varepsilon^{2h} + \varepsilon^{4h} + 4 \varepsilon^{3h} \cos(\theta)} +
2 \varepsilon^{h+1},$$
where $\theta$ is the angle pictured in Figure~\ref{F3}.

\begin{figure}[ht]\label{F3}
\psfrag{A}{$c({\bf x}, {T_{k}})$}
\psfrag{J}{$\varepsilon^{h}$}
\psfrag{I}{$\theta$}
\psfrag{H}{$c({\bf x}, {T_{i}})$}
\psfrag{E}{ }
\psfrag{F}{$\varepsilon^{h+1}$}
\psfrag{R}{$c({\bf x}, {T_{j}})$}
\psfrag{U}{ }
$$\includegraphics[width=10cm]{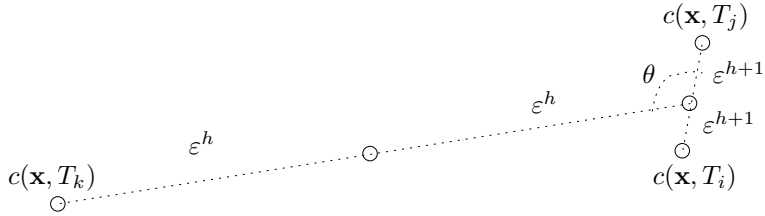}$$
\caption{ The geometry of $P_{T}$ for $T$ as in the Jacobi identity.}
\end{figure}

We argue by symmetry that the orientations of the\ three 
components of $\partial J$ gives rise to the Jacobi identity exactly.
Again referring to Figure~2, for each boundary component  we can define
an inward normal vector through having $c({\bf x}, {T_{i}})$ and $c({\bf x}, {T_{j}})$ move
radially outward, away from their center, and thus needing $c({\bf x}, {T_{k}})$ move
radially inward so that condition~(2.2) is satisfied.  This normal vector
is invariant under cyclic permutation of $T_{1}$, $T_{2}$ and $T_{3}$, as is
the definition of orientation for the $P_{T}$ for $T$ which appear in the
Jacobi identity.  Thus, these orientations will either all agree or all
disagree with a chosen
orientation of $J$, meaning in either case that the Jacobi identity holds.
\end{proof}

Finally, we allow for multiple planetary systems,  
freeing the points which do not move in the definition
of $P_{T}$, for example the points $x_{4}$ and $x_{5}$ in Figure~1.

\begin{definition}
\begin{itemize}
\item An $n$-forest is a collection of (by abuse) $S$-trees, with root vertices
at the points $(0, 0)$, $(1, 0)$, \ldots in the upper half-plane, where each
integer from $1$ to $n$ labels exactly one leaf.
\item If $F = \bigcup T_{i}$ is a forest, let $P_{F}$ be the submanifold defined 
by conditions~(\ref{p2}) and (\ref{p3}) of \refD{PT}, replacing condition~(\ref{p1}) with 
$c({\bf x}, {T_{i}}) = (i, 0, \ldots, 0)$.
\item Parameterize $P_{F}$ by a map of the same name from a product
over vertices of $F$ (ordered from left to right by the half-planar embedding
of $F$) of spheres, which when restricted to the coordinates labeled by $T_{i}$ is
a translation of $P_{T_{i}}$, namely $P_{T_{i}} + (i, 0, \ldots, 0)$.
\item By abuse let $F$ denote the homology class represented by $P_{F}$,
again using our fixed orientation of a sphere to orient the torus and thus its image
under $P_{F}$.
\end{itemize}
\end{definition}

We recover $P_{T}$ by letting $F$ be a forest which consists of $T$
and a collection of one-leaf trees.
We can summarize our results so far as follows.

\begin{definition}
Let $\poidn$ denote the quotient of the free module spanned by $n$-forests
by anti-symmetry and Jacobi identities
as in \refP{relations} along with the following:
\begin{enumerate}
\item[]  {\it{(commutativity) \, If $F_{1}$ and $F_{2}$ consist
of the same trees,  then $F_{1} = \sigma^{(d-1)} F_{2}$, where $\sigma$ is the
sign of the permutation which relates the ordering of the internal vertices of
the trees in $F_{1}$ with those of $F_{2}$.}}
\end{enumerate}
\end{definition}

\begin{theorem}
Sending a forest $F$ to the image of the fundamental class of $(S^{d-1})^{\times |F|}$
under $P_{F}$ gives a well-defined homomorphism from $\poidn$ to 
$H_{*}\left(\Co_{n}(\R^{d})\right)$.
\end{theorem}

Our main theorem will be that this map is an isomorphism (of {\em operads}).

\subsection{Canonical realization after compactification}

There are a number of choices we have made in our definition of $P_{T}$,
in particular the scale $\varepsilon$.  It is tempting to let $\varepsilon$
go to zero, which is indeed possible with some recent ``compactification technology.''
There is a canonical completion of these configuration spaces due to 
Fulton-MacPherson \cite{FuMa94} and Axelrod-Singer \cite{AxSi94},
which we denote $\Co_{n}[\R^{k}]$ in \cite{Sinh04}.  There we give the following
elementary definition.

\begin{definition}\label{D:compact}
\begin{itemize}
\item Define $\alpha_{ij} \colon \Co_n(\R^{d}) \to S^{d-1}$ as sending 
$(x_1, \cdots, x_n)$ to $\frac{x_j - x_i}{|x_j - x_i|}$.  
\item Let $I = [0,\infty]$, the one-point compactification of the nonnegative 
reals, and for $i,j,k$ distinct numbers between $1$ and $n$ 
let $s_{ijk} \colon \Co_n(\R^{d}) \to I = [0, \infty]$ 
be the map which sends $(x_1, \cdots, x_n)$ to  $\left(|x_i - x_j|/|x_i - x_k|\right)$. 
\item Let $\Co_{n}[\R^{d}]$ be the closer of the image of $\Co_{n}(\R^{d})$ 
in $$(\R^{d})^{\times n} \times \left(S^{d-1}\right)^{\times\binom{n}{2}} \times I^{\times \binom{n}{3}},$$
under the map which is the canonical inclusion on the first factor, the product over all $a_{ij}$ on the
second factor, and the product over all $s_{ijk}$ on the third factor.
\end{itemize}
\end{definition}

With such an elementary definition, the hard work is to establish basic properties.
This completion is functorial for embeddings and yields a manifold with corners with
$\Co_{n}(\R^{k})$ as its interior, and its strata are naturally indexed by trees (which are
not necessarily trivalent).  The points added by the boundary may be viewed as 
``degenerate configurations'' in which some number of points now coincide in the
large scale but have data which resolves them ``infinitesimally.''  See \cite{Sinh04} for a
more thorough treatment.

When the submanifold $P_{T}$ is included
in $\Co_{n}[\R^{k}]$,  we may send $\varepsilon$ in its definition to zero.  We leave this
as an exercise for the moment, since composing the map $P_{T}$ with the projections
$\alpha_{ij}$ will be a key calculation in \refT{agree}.
After this canonical homotopy of $P_{T}$ sending $\varepsilon$ to zero, the
resulting submanifold ends up being the stratum labeled by $T$ in the stratification
of $\Co_{n}[\R^{k}]$ as a manifold with corners.  
Thus these strata represent the homology classes we have
been constructing.   The Jacobi manifolds are also homotopic to strata 
which labeled by trees with one four-valent vertex, and their boundaries as strata
give rise to the Jacobi identity.  Indeed,
the manifold with boundary $\Co_{3}[\R^{k}]$ is diffeomorphic to 
the simplest Jacobi manifold.
Thus the compactification $\Co_{n}[\R^{k}]$  ``wears its homology 
on its strata.''  

\subsection{Historical notes}

To our knowledge, the  approach to homology through ``orbital systems''
first appeared implicitly in Fred Cohen's thesis work (\cite{CohF73}, in particular Section~12 of part III).  
This approach coincides with the ``twisted products''
constructions from Fadell and Husseini's book (\cite{FaHu01}, Chapter VI), but in their approach trees and forests 
do not explicitly appear.  The role of trees and forests, and the explicit connection with the theory of operads,
has been in the air since the ``renaissance'' of that theory, in particular \cite{GeJo94} which also emphasizes
the role of compactifications.  One of our aims is to
make this basic theory which is well-known to experts as accessible as possible.

\section{The cohomology ring}\label{S:cohomology}

In the previous section we constructed homology classes for the space
of Euclidean configurations
by mapping in fundamental classes of tori.  
In this section we pull back cohomology from tori.
If  homology classes in configuration spaces may be viewed as planetary
systems, cohomology classes may be view as recording planetary alignments.

\begin{definition}\label{D:aij}
Recall $\alpha_{ij} \colon \Co_n(\R^{d}) \to S^{d-1}$ as sending 
$(x_1, \cdots, x_n)$ to $\frac{x_j - x_i}{|x_j - x_i|}$.  
Let $\iota \in H^{d-1}(S^{d-1})$ denote the dual to the fundamental class, using
our fixed orientation. 
Let $a_{ij}$ denote $\alpha_{ij}^*(\iota)$. 
\end{definition}

The ring generated by these $a_{ij}$ can be represented graphically.

\begin{definition}
Consider graphs  with vertices labeled $1,\ldots,n$, 
with edges which are oriented and ordered.   Let  $\Gamma(n)$
denote the free module generated by such graphs, which is a ring by taking
the union of edges of two graphs in order to multiply them (using
the order of multiplication to define the ordering on the union of edges).  Map ${\Gamma}(n)$
to $H^*\left(\Co_n(\R^{d})\right)$ by sending a generator $\linep{i}{j}$ to $a_{ij}$.
\end{definition}

So for example the graph $ \linep{4}{2} \linep{1}{3}$, 
with $\linep{4}{2}$ first in the ordering of edges, is mapped to the product $a_{42} a_{13}$.
We will see that the map from $\Gamma(n)$ to $H^*\left(\Co_n(\R^{d})\right)$ is surjective.
As a base case, we show that after quotienting by the the relation $\linep{i}{j} =
(-1)^{d} \linep{j}{i}$, this map is an isomorphism in degree $d-1$.

\begin{definition}
Let $p_{i} : \Co_{n}(\R^{d}) \to \Co_{n-1}(\R^{d})$ be the projection map which sends
$(x_{1}, \ldots, x_{n})$ to $(x_{1}, \ldots, \hat{x_{i}}, \ldots, x_{n})$.
\end{definition}

\begin{lemma}
The projection map $p_{i}$ gives $\Co_{n}(\R^{d})$ the structure of a fiber bundle over
$\Co_{n-1}(\R^{d})$, with fiber given by $\R^{d}$ with $(n-1)$ points removed.
\end{lemma}

\begin{proof}
For simplicity let $i=n$.  Consider a neighborhood $U_{{\bf x}}$
of ${\bf x} = (x_{1}, \ldots, x_{n-1})$
of points $(y_{1}, \ldots, y_{n})$ where $d(y_{j}, x_{j}) < \epsilon$ for some fixed $\epsilon$
less than the minimum of the $\frac{1}{5}d(x_{j}, x_{k})$.  Construct a continuous family
of homeomorphsims $h_{{\bf y}}$ over $y \in U_{{\bf x}}$
between $\R^{d} - \{y_{1}, \ldots, y_{n-1}\}$ and $\R^{d} - \{B(x_{1}, \epsilon), \ldots,
B(x_{n-1}, \epsilon)\}$, which for good measure is the identity on 
$\R^{d} - \{B(x_{1}, 2\epsilon), \ldots,B(x_{n-1}, 2\epsilon)\}$.  This may be done, for
example, by ``straight-line retractions.''  
A trivialization of this fiber bundle, in other words a homeomorphism between 
$p_{n}^{-1}(U_{{\bf x}})$ and $U_{{\bf x}} \times \R^{d} - \{x_{1}, \ldots, x_{n-1}\}$ 
respecting $p_{n}$, is given by 
$$(y_{1}, \cdots, y_{n}) \mapsto (y_{1}, \cdots, y_{n-1}) \times h_{{\bf x}}^{-1} \circ h_{{\bf y}}(y_{n}).$$
\end{proof}

The space $\R^{d}$ with $(n-1)$ points removed retracts onto $\bigvee_{n-1}S^{d-1}$.
We assemble these projection maps into a tower of fibrations first studied by 
Fadell and Neuwirth \cite{FaNe62}, which is central
in the study of the topology of configuration spaces.
\begin{equation}\label{tower}
\xymatrix@=18pt{ 
   \bigvee_{n-1}S^{d-1}  \ar@<.5ex>[r]    & \Co_{n}(\R^{d}) \ar[d]^{p_{n}}\\
   \bigvee_{n-2}S^{d-1}  \ar[r]    & \Co_{n-1}(\R^{d}) \ar[d]^{p_{n-1}}\\
       & \vdots \ar[d]\\
   \bigvee_{2}S^{d-1}  \ar[r]    & \Co_{3}(\R^{d}) \ar[d]^{p_{n}}\\
      & \Co_{2}(\R^{d}) \simeq S^{d-1}
  }
\end{equation}
These fibrations split.  Choice of sections of $p_{i}$ include adding a new $i$th
point ``at infinity'' or somehow ``doubling'' the $i$th point.

\begin{proposition}
The first non-trivial homology group $H_{d-1}(\Co_{n}))$ is free of rank $\binom{n}{2}$.
\end{proposition}

\begin{proof}
We give a proof only for $n>2$.    We use the long exact sequences in homotopy of the fibrations 
constituting the tower of Equation~(\ref{tower}), which splits into short exact sequences because the
maps $p_{i}$ admit splittings.  We base all of our configuration spaces
at the configuration with $x_{i} = (i, 0, \ldots, 0)$.

From these short exact sequences, we deduce inductively 
that  $\Co_{i}(\R^{d})$ is $(d-2)$-connected,
and that $\pi_{d-1}\left(\Co_{i}(\R^{d})\right)$ is free.  Moreover, the rank of 
$\pi_{d-1}\left(\Co_{i}(\R^{d})\right)$
is that of $\pi_{d-1}\left(\Co_{i-1}(\R^{d})\right)$ plus $i-1$, or ultimately $1 + 2 + \cdots + (i-1)$,
which is $\binom{i}{2}$.  Finally, the Hurewicz theorem applies to give that the homology
is isomorphic to homotopy in this dimension.
\end{proof}

We now show that the homology and cohomology classes that we have constructed
so far generate this first non-trivial group.
Let $\la \, , \, \ra$ denote the standard pairing of cohomology and homology.

\begin{lemma}
Let $i<j$ and $k < \ell$.  The pairing $\la \linep{i}{j}, \tree{k}{\ell} \ra$ is
equal to one if $i = k$, $j = \ell$ and is zero otherwise.
\end{lemma}

\begin{proof}
To evaluate a cycle on $\linep{i}{j}$, it suffices by naturality to map that cycle
to $S^{d-1}$ and evaluate it on $\iota$.  Because 
$\tree{k}{\ell}$ is the natural image under $P_{\tree{k}{\ell}}$ 
of the fundamental class of $S^{d-1}$, it suffices to compute the degree 
of the composite $\alpha_{ij} \circ P_{\tree{k}{\ell}} : S^{d-1} \to S^{d-1}$.  
If $i = k$ and $j = \ell$, then this composite is the identity, which is degree one.
If we count the preimages of  the ``north pole'' in  $S^{d-1}$ to compute the degree, then
we are counting the number of configurations in $P_{\tree{k}{\ell}}$ for which 
$x_{j}$ is ``above'' $x_{i}$.  If either $i \neq k$ or $j \neq \ell$, then at no
point will $x_{j}$ be above $x_{i}$, since at least one of them will be
stationary, ``at infinity'', in every configuration parameterized by $P_{\tree{k}{\ell}}$.
\end{proof}

More generally, to evaluate a cohomology class of $\Co_{n}(\R^{d})$ 
represented by some graph, it suffices to count 
(with signs) the number of points in a cycle for which, for each edge $\linep{i}{j}$
in the graph
the point $x_{j}$ is ``above'' $x_{i}$. Since our cycles $P_{T}$ 
are planetary systems, the values of these cohomology classes on
them are counting planetary alignments.

Because $H_{d-1}\left(\Co_{n}(\R^{d})\right)$ is free of rank $\binom{n}{2}$, and this pairing
shows that the classes we have defined are linearly independent with the same
rank, we have the following.  

\begin{corollary}\label{C:evald-1}
A basis for $H^{d-1}(\Co_{n}(\R^{d}))$ is given by the classes $a_{ij}$.
The coefficients of a class expressed in this basis are given by evaluating
on the homology basis $(\tree{k}{\ell})$ with $k < \ell$.
\end{corollary}

\medskip

In general the map  from $\Gamma(n)$ to $H^*(\Co_n(\R^{d}))$ has relations.
If $G_{1}$ and $G_{2}$ differ by the reversal of
$k$ arrows and the reordering of edges as governed by a permutation $\sigma$,
then
\begin{equation*}\label{anti}
\text{\it (arrow reversing)} \qquad G_{1} - (-1)^{k(d-1)}({\rm{sign}}\ \sigma)^{d}G_{2} = 0
\end{equation*}
Also, because $\iota^{2} = 0$ in the cohomlogy
of the sphere, any graph with more than one edge between two vertices will
map to zero.  There is a more subtle relation, which is in some sense
dual to the Jacobi identity.

\begin{theorem}\label{T:arnold}
The following relation holds in the image of $\Gamma(n)$ in $H^*(\Co_n(\R^{d}))$:
$$
\text{(Arnold)}\qquad  \qquad
\begin{xy}                           
  (0,-2)*+UR{\scriptstyle j}="a",    
  (3,3)*+UR{\scriptstyle k}="b",   
  (6,-2)*+UR{\scriptstyle \ell}="c",   
  "a";"b"**\dir{-}?>*\dir{>},         
  "b";"c"**\dir{-}?>*\dir{>},         
  (3,-5),{\ar@{. }@(l,l)(3,6)},
  ?!{"a";"a"+/va(210)/}="a1",
  ?!{"a";"a"+/va(240)/}="a2",
  ?!{"a";"a"+/va(270)/}="a3",
  ?!{"b";"b"+/va(120)/}="b1",
  "a";"a1"**\dir{-},  "a";"a2"**\dir{-},  "a";"a3"**\dir{-},
  "b";"b1"**\dir{-}, "b";(3,6)**\dir{-},
  (3,-5),{\ar@{. }@(r,r)(3,6)},
  ?!{"c";"c"+/va(-90)/}="c1",
  ?!{"c";"c"+/va(-60)/}="c2",
  ?!{"c";"c"+/va(-30)/}="c3",
  ?!{"b";"b"+/va(60)/}="b3",
  "c";"c1"**\dir{-},  "c";"c2"**\dir{-},  "c";"c3"**\dir{-},
  "b";"b3"**\dir{-}, 
\end{xy}\ + \                             
\begin{xy}                           
  (0,-2)*+UR{\scriptstyle j}="a",    
  (3,3)*+UR{\scriptstyle k}="b",   
  (6,-2)*+UR{\scriptstyle \ell}="c",    
  "b";"c"**\dir{-}?>*\dir{>},         
  "c";"a"**\dir{-}?>*\dir{>},          
  (3,-5),{\ar@{. }@(l,l)(3,6)},
  ?!{"a";"a"+/va(210)/}="a1",
  ?!{"a";"a"+/va(240)/}="a2",
  ?!{"a";"a"+/va(270)/}="a3",
  ?!{"b";"b"+/va(120)/}="b1",
  "a";"a1"**\dir{-},  "a";"a2"**\dir{-},  "a";"a3"**\dir{-},
  "b";"b1"**\dir{-}, "b";(3,6)**\dir{-},
  (3,-5),{\ar@{. }@(r,r)(3,6)},
  ?!{"c";"c"+/va(-90)/}="c1",
  ?!{"c";"c"+/va(-60)/}="c2",
  ?!{"c";"c"+/va(-30)/}="c3",
  ?!{"b";"b"+/va(60)/}="b3",
  "c";"c1"**\dir{-},  "c";"c2"**\dir{-},  "c";"c3"**\dir{-},
  "b";"b3"**\dir{-}, 
\end{xy}\ + \                              
\begin{xy}                           
  (0,-2)*+UR{\scriptstyle j}="a",    
  (3,3)*+UR{\scriptstyle k}="b",   
  (6,-2)*+UR{\scriptstyle \ell}="c",    
  "a";"b"**\dir{-}?>*\dir{>},         
  "c";"a"**\dir{-}?>*\dir{>},          
  (3,-5),{\ar@{. }@(l,l)(3,6)},
  ?!{"a";"a"+/va(210)/}="a1",
  ?!{"a";"a"+/va(240)/}="a2",
  ?!{"a";"a"+/va(270)/}="a3",
  ?!{"b";"b"+/va(120)/}="b1",
  "a";"a1"**\dir{-},  "a";"a2"**\dir{-},  "a";"a3"**\dir{-},
  "b";"b1"**\dir{-}, "b";(3,6)**\dir{-},
  (3,-5),{\ar@{. }@(r,r)(3,6)},
  ?!{"c";"c"+/va(-90)/}="c1",
  ?!{"c";"c"+/va(-60)/}="c2",
  ?!{"c";"c"+/va(-30)/}="c3",
  ?!{"b";"b"+/va(60)/}="b3",
  "c";"c1"**\dir{-},  "c";"c2"**\dir{-},  "c";"c3"**\dir{-},
  "b";"b3"**\dir{-}, 
\end{xy}\ =\ 0,
$$
where ${j}$, ${k}$, and ${\ell}$ stand for 
vertices in the graph which could possibly have other connections to
other parts of the graph (indicated by the ends of edges abutting 
$j$, $k$, and $\ell$)
which are not modified in these operations.  
\end{theorem}

\begin{proof}
Using the ring structure, it suffices to consider when 
there there are no edges incident on $j$, $k$
and $\ell$, other than the two edges involved in the identity.  
In this case the cohomology classes are all pulled 
back from $H^{2(d-1)}\left(\Co_{3}(\R^{d})\right)$ via a map 
which forgets all $x_{i}$ except $x_{j}$, $x_{k}$ and $x_{\ell}$.
So we may assume that $n=3$ and $\{j, k, \ell\} = \{1, 2, 3\}$.

Our proof uses elementary intersection theory to compute some cup products.
Since $\Co_{3}(\R^{k})$ is a manifold, its cohomology is Lefshetz dual to 
its locally finite homology.  
Consider the submanifold of 
$(x_{1}, x_{2}, x_{3}) \in \Co_{3}(\R^{k})$   such that
$x_{1}$, $x_{2}$, and $x_{3}$ are collinear.
This submanifold has three components.  Let $\col_{i}$ denote the component in 
which $x_{i}$ is in the middle.  
Since $\col_{i}$ is a properly embedded submanifold of codimension $d-1$, 
once oriented it represents a locally finite homology class, which through Lefshetz duality 
gives rise to a class in $H^{d-1}(\Co_{3}(\R^{d}))$.  By \refC{evald-1}, 
this class is determined by its value on the classes $\tree{j}{k}$, which in the context
of Lefshetz duality means intersecting $\col_{i}$ with 
various $P_T$.

These intersections can be understood directly.  The submanifold 
$\col_{i}$ can only intersect $P_{\tree{i}{j}}$ or $P_{\tree{i}{k}}$, since otherwise $x_{i}$
would be ``at infinity'', and thus
could not be in the middle of a collinearity.  Moreover, $P_{\tree{i}{j}}$ and $P_{\tree{i}{k}}$ each
intersect $\col_{i}$ exactly once, namely when $x_{i}$ is a negative multiple of
$x_{j}$ (respectively, $x_{k}$).  For purposes of our computation we only need
that these intersections differ in sign by $-1$, coming from orientation reversing
of the line on which the three points lie.  We deduce that the
cohomology class represented by $\col_{i}$ is $\pm(a_{ij} - a_{ik})$.

Under Lefshetz duality, cup products are computed by
(transversal) intersections.  Since $\col_{1}$ and $\col_{2}$ are disjoint, the
cohomology classes which they represent cup to zero.  We have
\begin{align*}
0 &= \pm (a_{12} - a_{13})(a_{23} - a_{21}) \\
&= a_{12}a_{23} - a_{12}a_{21} - a_{13}a_{23} + a_{13}a_{21} \\
&= a_{12}a_{23} + 0 - (-1)^{d + (d+1)}a_{23}a_{31} +(-1)^{2d} a_{31}a_{12} \\
&=  a_{12}a_{23} + a_{23}a_{31} + a_{31}a_{12}.
\end{align*}
When we translate back to the graphical language, this is exactly the Arnold
identity.
\end{proof}

This new proof through the disjointness of collinearity submanifolds is using
a fundamental geometric observation as the basis for a cohomology ring computation,
akin to seeing the cohomology ring of projective spaces through the 
intersections of linear subspaces.  
Using cochains defined through collinearities works better for
this calculation than our
original cochains representing the classes $a_{ij}$, for which we have that
the quadratic polynomial in the Arnold identity is exact but not identically
zero (unless $d=2$, where Kontsevich observed the vanishing using the
differential forms $\dlog(x_{i} - x_{j})$.)  The collinearity cochains are also invaraint
under the action of the full group of affine transformations in $\R^{d}$.

We will see that there are no further relations among these graph classes
in $H^{*}\left(\Co_{n}(\R^{d})\right)$, so that the image of $\Gamma(n)$ 
will be precisely the following module.

\begin{definition}
Let $\iopd(n)$ denote the quotient of $\Gamma(n)$ by the
arrow-reversing relation and the Arnold identity.
\end{definition}

Instead of starting with $\Gamma(n)$, we can restrict to acyclic graphs.

\begin{proposition}
Any element of $\iopd(n)$ represented by a graph which has a cycle is zero.
\end{proposition}

\begin{proof}
We may use the Arnold identity inductively to reduce to graphs with shorter cycles.
But graphs with cycles of length two, that is which have more than one edge between
two vertices, are zero.
\end{proof}

\subsection{Historical notes}
Analyzing projection maps and in particular assembling them into a tower has been a central
tool in this area since its introduction by Fadell and Neuwirth \cite{FaNe62}.
The calculation of cohomology of configurations in the plane is, famously, due to Arnold \cite{Arno69}.
It was generalized in higher dimensions by Cohen \cite{CLM76}, to the complements of other collections
of subspaces defined by linear equations by Orlik and Solomon \cite{OrSo80}, and to a myriad of
other contexts.   Note that what we call the Arnold identity, along with the fact that generators square to zero, 
are called the cohomological Yang-Baxter relations in Chapter V of \cite{FaHu01}, which gives a complete
account of the spectral sequence approach to calculation.  The graphical notation for this cohomology has
been useful for a wide range of recent work, for example \cite{KoSo00}.

\section{The homology-cohomology pairing}\label{S:pairing}

Building on the combinatorics which arose in the last two sections, we develop 
a pairing between graphs and trees which coincides with the
evaluation of cohomology on homology of $\Co_{n}(\R^{d})$.

\begin{definition}\label{D:pairing}
Given an $n$-graph $G$ and an $n$-tree $T$, define the map 
$$\beta_{G,T}:\bigl\{\text{edges of } G\bigr\} \longrightarrow 
\bigl\{\text{internal vertices of } T\bigr\}$$ 
by sending an edge $\linep{i}{j}$ in $G$ to the vertex at the nadir of 
the shortest path in $T$ between the leafs with labels $i$ and  $j$. 
Define the mod-$2$ configuration pairing of $n$-graphs 
$G$ and $n$-trees $T$ by 
$$\bigl\langle G,\, T\bigr\rangle = 
  \begin{cases} 1                & \text{if $\beta$ is a bijection} \\
   0  & \text{otherwise}.
\end{cases}$$
Define the dimension-$d$ configuration pairing by, in the first case above,
introducing the sign of the permutation relating the orderings of the edges
of $G$ and internal vertices of $T$ given in \refD{trees}
when $d$ is even, or $(-1)^{k}$ where $k$
is the number of edges $\linep{i}{j}$ of $G$ for which leaf $i$
is to the right of leaf $j$ under the planar embedding of $T$ when $d$ is odd.

 This definition extends to give a pairing
 between (possibly disconnected) $n$-graphs $G$  and $n$-forests $F$,
which is zero if an edge of $G$ has endpoints which label leaves in two different 
components of $F$ (so that $\beta$ is not defined).
\end{definition}

\begin{figure}[ht]
$$\begin{aligned}\begin{xy} 
  (0,-2.5)*+UR{\scriptstyle 1}="a", 
  (3.75,3.75)*+UR{\scriptstyle 2}="b", 
  (7.5,-2.5)*+UR{\scriptstyle 3}="c", 
  "a";"b"**\dir{-}?>*\dir{>} ?(.4)*!RD{\scriptstyle e_1}, 
  "b";"c"**\dir{-}?>*\dir{>} ?(.5)*!LD{\scriptstyle e_2} 
 \end{xy}\end{aligned}\  \longmapsto \ 
 \begin{xy}
   (2,1.5); (4,3.5)**\dir{-}, 
   (0,3.5); (4,-.5)**\dir{-} ?(.5)*!RU{\scriptstyle \beta(e_1)}; 
   (4,-.5); (8,3.5)**\dir{-} ?(.2)*!LU{\scriptstyle \beta(e_2)}, 
   (4,-.5); (4,-2.5)**\dir{-}, 
   (0 ,4.7)*{\scriptstyle 2}, 
   (4,4.7)*{\scriptstyle 1}, 
   (8,4.7)*{\scriptstyle 3}, 
   (2,1.5)*{\scriptstyle \bullet},
   (4,-.5)*{\scriptstyle \bullet},
 \end{xy}  \qquad \qquad \qquad 
 \begin{aligned}\begin{xy} 
  (0,-2.5)*+UR{\scriptstyle 1}="a", 
  (3.75,3.75)*+UR{\scriptstyle 2}="b", 
  (7.5,-2.5)*+UR{\scriptstyle 3}="c", 
  "a";"b"**\dir{-}?>*\dir{>} ?(.4)*!RD{\scriptstyle e_1}, 
  "b";"c"**\dir{-}?>*\dir{>} ?(.5)*!LD{\scriptstyle e_2} 
 \end{xy}\end{aligned}\  \longmapsto \  
 \begin{xy}
   (2,1.5); (4,3.5)**\dir{-}, 
   (0,3.5); (4,-.5)**\dir{-} ?(.8)*!RU{\scriptstyle \beta(e_1)}; 
   (4,-.5); (8,3.5)**\dir{-} ?(.2)*!LU{\scriptstyle \beta(e_2)}, 
   (4,-.5); (4,-2.5)**\dir{-}, 
   (0 ,4.7)*{\scriptstyle 1}, 
   (4,4.7)*{\scriptstyle 3}, 
   (8,4.7)*{\scriptstyle 2}, 
   (4,-.5)*{\scriptstyle \bullet}, 
 \end{xy} $$
 \caption{The map $\beta_{G,T}$ for two different trees $T$.  In the first
 case the configuration 
 pairing is $-1$ if $d$ is odd or $1$ if $d$ is even, and in the second case it is zero.}
\end{figure}
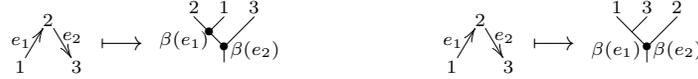

\begin{theorem}\label{T:agree}
The homology-cohomology pairing for $\Co_{n}(\R^{d})$ agrees with the configuration pairing.
That is, if we let $\la -, - \ra_{c}$ denote the combinatorially-defined configuration pairing,
and let $\la -, - \ra_{H}$ denote the homology-cohomology pairing for
$\Co_{n}(\R^{d})$, then $\la G, F \ra_{c} = \la G, F \ra_{H}$.  
\end{theorem}

Here we have continued the abuse of letting $G$ and $F$ denote both
graphs and trees and their images in cohomology and homology of $\Co_{n}(\R^{d})$.

\begin{proof}
For the homology pairing, we must evaluate a product of 
classes $a_{ij}$ on a submanifold $P_{F}$,
which by naturality of cap products is equal to computing the degree of the composite 
$$ \pi_{G} \circ P_{F} : \prod_{v \in F} S^{d-1} \overset{P_{F}}{\longrightarrow}
 \Co_{n}(\R^{d}) \overset{\pi_{G}}{\longrightarrow} \prod_{e \in G}S^{d-1}.$$
Here $\pi_{G}$ is the product over edges of $G$ which associates to
$e = \linep{i}{j}$ a factor of  $\pi_{ij}$.
By Definitions~\ref{D:param} and \ref{D:aij}, 
this composite sends 
$(u_{v_{1}}, \ldots, u_{v_{|F|}})$ to $(\theta_{e_{1}}, \ldots, \theta_{e_{|G|}})$,
where for  $e = \linep{i}{j}$, $\theta_{e}$ is the unit vector in the direction of
$$\left( (n,0,\cdots,0) + \sum_{v_{k}<{\rm leaf} \; i} \pm \varepsilon^{h_{k}} u_{v_{k}} \right) -
\left( (m, 0, \dots, 0) +  \sum_{w_{\ell}<{\rm leaf} \; j} \pm \varepsilon^{h_{\ell}} u_{w_{\ell}} \right).$$
Here $v_{k}$ is the vertex of height $k$ under leaf $i$, which is in the $n$th component
of the forest $F$,  and similarly $w_{\ell}$
is the vertex of height $\ell$ under leaf $j$, which is in the $m$th component of $F$.
If leaves $i$ and $j$ are in the same component, common terms associated to 
vertices under both leaf $i$
and leaf $j$ cancel, leaving $\theta_{e}$ as the unit vector in the direction of
$$
\varepsilon^{h_{n}} \left(\pm 2 u_{v_{n}} \pm \varepsilon (u_{v_{n+1}} - u_{w_{n+1}}) +
     \varepsilon^{2}     \cdots \right),$$
where $v_{n}$ is the highest vertex under both leaf $i$ and $j$ (if it exists), which
is also the nadir of the path between $i$ and $j$.  

Consider the homotopy of $P_{F}$, and thus this composite, in which  
$\varepsilon$  approaches zero.   Through this homotopy, 
$\theta_{e}$ approaches either $(\pm 1,0, \cdots, 0)$ if leaves 
$i$ and $j$ are in different components, or otherwise $\sigma_{e} u_{v_{n}}$.
From \refD{param} we see
that $\sigma_{e}$ is $1$ if leaf $i$ is to the left of $j$ or $-1$ if it is to the right.
Therefore, if there is some edge $\linep{i}{j}$ in $G$ with leaves $i$
and $j$ in different components of $F$, then $P_{F}$ is homotopic to the map between
tori with at least one factor the constant map of $(\pm 1, \cdots, 0)$, and thus it is of degree zero.
Otherwise, $P_{F}$ is homotopic to the map which sends
$(u_{v_{1}}, \ldots, u_{v_{|F|}}) \in \prod_{|F|} S^{d-1}$ to 
$$\left(\sigma_{e_{1}}u_{\beta_{G,F}(e_{1})}, \ldots, 
  \sigma_{e_{|G|}}u_{\beta_{G,F}(e_{|G|})} \right),$$
whose degree agrees with the definition of the configuration pairing through 
$\beta_{G,F}$.
\end{proof}

Because the homology classes of $\Co_{n}(\R^{d})$ represented by forests satisfy the 
anti-symmetry and Jacobi identities of \refP{relations}, and the cohomology 
classes represented by graphs  
satisfy arrow-reversing and the Arnold identities, we have geometrically 
established the following fact, which is established combinatorially in \cite{Sinh06.2}.

\begin{corollary}\label{L:1}
The dimension-$d$ configuration pairing passes to a well-defined pairing
between $\poidn$ and $\iopd(n)$.
\end{corollary}

We now outline the purely algebraic argument, given in \cite{Sinh06.2}, that
this pairing between $\poidn$ and $\iopd(n)$ is perfect.
We only give hints, leaving some fun for the reader.

\begin{lemma}\label{L:3}
The module $\poidn$ is spanned by $n$-forests in which all trees are {\em tall}
(that is, the distance between the leaf with the minimal label and the root is maximal,
and that leaf is leftmost in the planar ordering).
The module $\iopd(n)$ is spanned by $n$-graphs whose components
are {\em long} (that is, each component is a linear graph, with one endpoint labeled by
the minimal label; edges are ordered consecutively and 
oriented away from the minimal label).
\end{lemma}

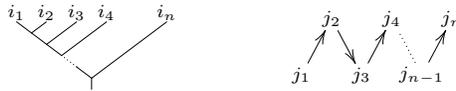
\begin{figure}[ht]\label{tallnlong}
$$\begin{aligned}\begin{xy}
   (2,1.5); (4,3)**\dir{-},	
   (4,0); (8,3)**\dir{-};
   (0,3); (6,-1.5)**\dir{-};
   (12,3)**\dir{-},
   (6,-1.5);(8,-3)**\dir{.};
   (10,-4.5)**\dir{-};
   (20,3)**\dir{-},
   (10,-4.5);(10,-6)**\dir{-},
   (0,4.2)*{\scriptstyle{i_{1}}},
   (4,4.2)*{\scriptstyle{i_2}},
   (8,4.2)*{\scriptstyle{i_3}},
   (12,4.2)*{\scriptstyle{i_4}},
   (20,4.2)*{\scriptstyle{i_n}}
 \end{xy}\end{aligned} \qquad \qquad
 \begin{aligned}\begin{xy}
   (0,-3)*+UR{\scriptstyle j_{1}}="1",
   (4,4.5)*+UR{\scriptstyle j_2}="2",
   (8,-3)*+UR{\scriptstyle j_3}="3",
   (12,4.5)*+UR{\scriptstyle j_4}="4",
   (16,-3)*+UR{\scriptstyle j_{n-1}}="5",
   (20,4.5)*+UR{\scriptstyle j_n}="6",
   "1";"2"**\dir{-}?>*\dir{>},
   "2";"3"**\dir{-}?>*\dir{>},
   "3";"4"**\dir{-}?>*\dir{>},
   "4";"5"**\dir{.},
   "5";"6"**\dir{-}?>*\dir{>},
 \end{xy}\end{aligned}$$
 \caption{A tall tree and a long graph.  Here $i_{1}$ and $j_{1}$ are minimal among
 indices in the tree and graph respectively.}
\end{figure}

\begin{proof}[Sketch of proof] 
For the forests, use the Jacobi identity inductively to increase the distance from the minimally
labeled leaf to the root.  For the graphs, use the Arnold identity to reduce the number
of edges incident on a given vertex.
\end{proof}

The sets of  tall forests and long graphs are in one-to-one correspondence with 
partitions of $\n$, where
both the subsets and the constituent elements of each subset are ordered.

\begin{lemma}\label{L:4}
The degree-$d$ configuration pairing of a tall forest and a long graph 
is equal to one if their associated ordered
partitions agree, and is zero otherwise.
\end{lemma}

\begin{proof}[Sketch of proof] 
By definition, the underlying unordered partitions must agree in order for the pairing
to be non-zero.  When looking at the configuration pairing between a single tall tree $T$
and long graph $G$ which
share labels, look at the first place where their orderings differ to see how
$\beta_{G,T}$ fails to be a bijection.
\end{proof}

Thus, on tall forests and long graphs, the configuration pairing is essentially a Kronecker
pairing, showing that these spanning sets are linearly independent.

\begin{corollary}\label{C:5}
Tall forests form a basis of $\poidn$.  Long graphs form a basis of $\iopd(n)$.
Both $\poidn$ and $\iopd(n)$ are torsion-free.
\end{corollary}

Because tall forests and long graphs form bases, and the dimension-$d$
configuration pairing is a Kronecker pairing on them, we deduce the main algebraic
result.

\begin{theorem}\label{T:perfect}
The dimension-$d$ configuration pairing between $\poidn$ and $\iopd(n)$ 
is perfect.
\end{theorem}

And because the configuration pairing agrees with the homology pairing for 
$\Co_{n}(\R^{d})$ on the classes constructed by graphs and forests, we have the following.

\begin{corollary}\label{C:inj}
The homomorphisms from $\poidn$ to $H_{*}\left(\Co_{n}(\R^{d})\right)$ and 
$\iopd(n)$ to $H_{*}\left(\Co_{n}(\R^{d})\right)$ are injective.
\end{corollary}

We can now establish the first part of the main result of this paper.

\begin{theorem}\label{T:main1}
The maps from $\poidn$ to $H_{*}(\Co_{n}(\R^{d}))$ and 
$\iopd(n)$ to $H^{}{*}(\Co_{n}(\R^{d}))$ are isomorphisms.
\end{theorem}

\begin{proof}
We make light
use of the  Leray-Serre spectral sequence.  If in a  fibration 
$F \to E \to B$, we have that  $B$ is simply connected and either $F$ or $B$ has torsion-free
homology, then this spectral sequence says that $H^{*}(F) \otimes H^{*}(B)$
serves as an ``upper bound'' for $H^{*}(E)$.  That is,
if we let $P(X)(t) = \sum \; ({\rm rank} \; H_{i}(X)) t^{i}$, 
the  Poincar\'e polynomial of $X$, 
then $P(E) \leq P(F) \cdot P(B)$, where by $\leq$ we mean that this inequality
holds for all coefficients 
of $t^{i}$.  Moreover, if the homology of $F$ and $B$ are free, equality 
is achieved only when that of $E$ is free.

Recall the Fadell-Neuwirth tower of fibrations from Equation~(\ref{tower}).  If $d>2$ then the long exact sequence
for homotopy groups for the fibration $\bigvee_{i-1} S^{d-1} \to \Co_{i}(\R^{d})  \to \Co_{i-1}(\R^{d})$ can be used to inductively establish that these configuration spaces are simply connected.
The Leray-Serre spectral sequence upper bound then yields the
inequality $P( \Co_{i}(\R^{d})) \leq P( \Co_{i-1}(\R^{d}) \cdot (1 + (i-1) t^{d-1})$.
For $d=2$ these spaces are not simply connected - in fact their fundamental groups are pure 
braid groups almost by definition, and in fact they are classifying spaces for pure braid groups \cite{FoNe62}  -
but Fred Cohen does the extra work at the beginning of Part III of \cite{CLM76} necessary to show that this
upper bound still holds.

Inductively we have 
$$P\left( \Co_{n}(\R^{d})\right) \leq  \prod_{i = 1}^{n-1} (1 + i t^{d-1}).$$ 
We claim that this upper bound is sharp.   Let $Q_{n}$ be  polynomial
defined as $Q_{n}(t) = q_{i} t^{i(d-1)}$, where $q_{i}$ is the rank of
the submodule of $\iopd(n)$ with $i$ edges.  
By \refC{inj}, $Q_{n}$ is a lower bound for 
$P( \Co_{n}(\R^{d}))$, which we compute inductively.  
The set of long graphs in $\iopd(i)$
maps to those of $\iopd(i-1)$ by taking a long graph, removing the
vertex labeled $n$ and any edges connected to it, and then reconnecting the two
adjacent vertices with a new edge if necessary.  Given a long graph $G$ in 
$\iopd(i-1)$ there is exactly one long graph in $\iopd(i)$
with the same number of edges which maps to it, namely the one in 
which vertex $n$ is added but not connected to an edge.  Moreover,
there are $i-1$ long graphs in $\iopd(i)$ with one more edge
which map to it, since one can choose which of the $i-1$ vertices
in $G$ would have an edge connect to (as opposed to from) the $i$th vertex.
We deduce that $Q_{i} = Q_{i-1} \cdot  (1 + (i-1) t^{d-1})$.

Thus, the lower bound for $H^{*}(\Co_{n}(\R^{d}))$ given by the
submodule $\iopd(n)$ matches the upper bound given by
the Leray-Serre spectral sequences for the tower of fibrations of 
Equation~(\ref{tower}).  We may inductively deduce
that $H^{*}(\Co_{n}(\R^{d}))$ is free, isomorphic to $\iopd(n)$.
By the Universal Coefficient Theorem and Theorems~\ref{T:agree} and 
\ref{T:perfect}, we have that $H_{*}(\Co_{n}(\R^{d}))$ is isomorphic
to $\poidn$.
\end{proof}

One can also obtain upper bounds by calculations in the cellular homology of the
one-point compactifications of these spaces, which is then dual to their chohomology,
using a cell structure first developed for $d=2$ by Fox and Neuwirth \cite{FoNe62}.

Using the Leray-Serre spectral sequence in full force can lead to some
of these results more quickly.  For formal reasons, the spectral sequence
for each Fadell-Neuwirth fibration collapses, showing
immediately that the upper bound on cohomology groups 
given by each spectral sequence is sharp.  One
can also use a symmetry argument to deduce the Arnold identity and thus
determine the cohomology ring structure.  Indeed, these fiber sequences are
nice first examples to work with, since even though 
the group structure mimics that of a trivial
(product) fiber sequence, the cohomology ring of $\Co_{n}(\R^{d})$ differs
greatly from that of $\prod_{i=1}^{n-1} \left( \bigvee_{i} S^{d-1} \right)$.

In our approach, we not only have an understanding of the  homology groups 
of $\Co_{n}(\R^{d})$ and the cohomology ring up to isomorphism, but 
we also have canonical spanning sets and an explicit understanding
of the pairing between them, which enables hands-on calculations.

\subsection{Historical notes}
The pairing between graphs and trees we develop (for $d$ odd) was first noticed by Melan{\c{c}}on
and Reutenauer in a combinatorial study of free Lie superalgebras \cite{MeRe96}.  It was
independently identified as the pairing between canonical spanning sets for homology and cohomology
of configuration spaces by Paolo Salvatore, Victor Tourtchine  \cite{Tour04}, and the author, all within the context of
studying spaces of knots.  The present paper gives the first full account of this connection, to our knowledge.  
The pairing is useful in related areas of algebra and topology, such as the study of Hopf invariants
\cite{SiWa08}.

\section{Operads}\label{S:operads}

An operad encodes multiplication.  Roughly speaking, an operad contains
information needed to multiply in an algebra over that operad.
For example, in multiplying matrices one must supply an ordering of the 
matrices to be multiplied, while in multiplying real numbers no such 
ordering is needed.  To Lie multiply (that is, take commutators of)
some matrices, one must not only order but parenthesize them.

The many definitions of an abstract operad 
are necessarily complicated.  Even the elegant ``an operad is a monoid in the category
of symmetric sequences,'' requires knowing what a symmetric sequence is and
then doing some work to relate that definition to standard examples.  
Thorough introductions
to the theory of operads are given elsewhere in this volume.  
We prefer to be self-contained and to work with operads through trees, so we give
our own development here.  
We start with examples before giving the definition.  For now
we work with the intuitive definition that an operad $\Op$
in a symmetric monoidal category $\mathcal{C}$ 
is a sequence of objects indexed by natural
numbers so that the $n$th object $\Op(n)$ 
parameterizes ways in which $n$ elements
of some kind of algebra (that is, an algebra over that operad) can be multiplied.

\begin{examples}\label{Ex:operads}
\begin{enumerate}
\item In any unital symmetric monoidal category $\mC$, the commutative operad 
$\com$ has $\com(n) = \one_{\mC}$, since there is only one way to multiply
$n$ things commutatively.  (For vector spaces, $\one_{\mC}$ is the ground
field $\ka$; for spaces, $\one_{\mC}$ is a point.)  
\item In spaces, the associative operad $\ass$ has $\ass(n) = \Sigma_{n}$,
the finite set of orderings of $n$ points, since the product of $n$ things is
determined by their order if multiplication is associative.  In vector spaces,
$\ass(n) = \ka[\Sigma_{n}]$.
\item In vector spaces, the Lie operad $\lie$ has $\lie(n)$ spanned by $n$-trees
modulo the anti-symmetry and Jacobi identities of \refP{relations} (with $d$ even).
In a Lie algebra one must parenthesize elements to multiply them, and 
our $n$-trees encode parenthesizations.  The anti-symmetry
and Jacobi identities are always respected in a Lie algebra, so they appear in the 
definition of this operad.
\item  In spaces, the little $d$-disks operad has $\dis^{d}(n)$ the space of $n$ 
little disks in $D^d$.  Explicitly, $\dis^{d}(n)$ is the subspace of $\Co_n({\rm Int}D^d) \times (0,1]^n$ of 
$(x_i) \times (r_i)$ such that the balls $B(x_i, r_i)$ are contained in the interior of 
$D^d$ and have disjoint interiors.   This space parameterizes some ways in which 
maps from $S^{d}$, or rather $D^{d}$ with its boundary going to the basepoint,
can be multiplied.  Given $n$ maps $f_{i}: (D^{d}, \partial) \to (X, *)$ 
and an element of this space we may multiply the $f_{i}$ by applying them on the
corresponding little disk, sending points outside any little disk to the basepoint of $X$.
See Figure~7.
\item We've already seen the degree-$d$
Poisson operad in vector spaces, which has $n$th
entry $\poidn$.  In a Poisson algebra one can both Lie multiply, represented
by trees, and then multiply those results together, represented by placing those
trees together in a forest.  One can of course also multiply first and then Lie multiply.
By the Leibniz rule, the bracket is a derivation with 
respect to the multiplication, so we may always reduce to expressions which
correspond to forests.  For example (using bracket notation rather than forests and
assuming $d$ is even),
\begin{align*}
[x_{1}, [x_{2}, x_{3} \cdot x_{4}]] &= [x_{1}, [x_{2}, x_{3}] \cdot x_{4}]] + 
[x_{1}, x_{3} \cdot [x_{2}, x_{4}]] \\
&= [x_{1}, [x_{2}, x_{3}]] \cdot x_{4}  
+ [x_{1}, x_{4}] \cdot [x_{2}, x_{3}] + [x_{1}, x_{3}] \cdot [x_{2}, x_{4}] + [x_{1}, [x_{2}, x_{4}]] \cdot x_{3}.
\end{align*}
\item For any $X$ in $\mC$, the endomorphism operad of $X$ has $\ndo_{X}(n) =
{\rm Hom}_{\mC}(X^{\odot n}, X)$.  The endomorphism operad is often too large to understand
explicitly, much as groups of homeomorphisms are.  But finding an interesting sub-operad of the
endomorphism operad will endow $X$ with a multiplication, much as finding a sub-group
of the self-equivalences of $X$ gives a group action.
\end{enumerate}
\end{examples}

As McClure and Smith
point out  in \cite{McSm04.2}, the axiomatic definitions of operads follow nicely
from reflecting on what they do, just as the axioms of a group all follow
from the notion that a group encodes symmetries through a group action. 
For example, if one has a rule for multiplying four inputs and two inputs, then
one can make a rule for multiplying five inputs:   first multiply two of them,
then take that result as an input with the remaining three original inputs in order
to apply the rule for multiplying four inputs.  Thus, for any operad there
are maps $\Op(4) \odot \Op(2) \to \Op(5)$.

\begin{center}
\begin{figure}[ht]\label{F:basicmor}
\begin{center}
\psfrag{C1}{}
\psfrag{C2}{}
$$\includegraphics[width=11cm]{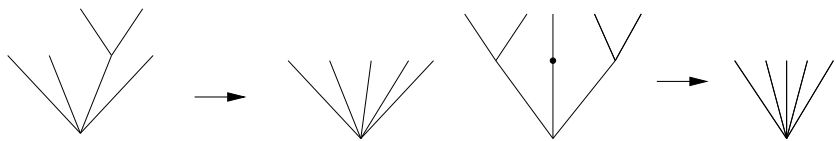}$$
\caption{Two morphisms in $\Upsilon$ which give rise to standard operad structure maps.
The first corresponds to a $\circ_i$ operation, the second to one of May's structure maps.}
\end{center}
\end{figure}
\end{center}

We prefer to state
the definition using trees, and, because operads require maps which
commute, the language of categories and functors.  In this language,
what we just said about combining two- and four-input multiplications to 
make a five-input multiplication is encoded in the first part of Figure~\ref{F:basicmor}.

\begin{definition}
\begin{itemize}
\item An o-tree is a finite connected acyclic graph with a 
distinguished vertex called the
root.  Univalent vertices of an o-tree (not counting the root, if it is 
univalent) are called leaves. 

\item At each vertex, the edge which is closer to the root
is called the output edge.  The edges which are further from
the root are called input edges and are labeled from
$1$ to $n$.

\item At each edge, the vertex of
an edge which is further from the root is called its input vertex, and
the vertex closer to the root is called its output vertex.  We say
that one vertex or edge lies over another if the latter is in the
path to the root of the former.  A non-root edge is called redundant if its 
initial vertex is bivalent.

\item Given an o-tree $\tau$ and an edge $e$, the
contraction of $\tau$ by $e$ is the tree $\tau'$ obtained by
identifying the input vertex of $e$ with its
output vertex,  and removing $e$ from the set of
edges.  If the label of $e$ was $i$, the labels of the 
$k$ edges which were immediately over $e$ will be increased
by $i-1$, and the labels of the edges which shared the output
vertex of $e$ with labels greater than $i$ will be increased by $k-1$.

\item Let $\Upsilon$ denote the category whose objects are o-trees, 
and whose morphisms are generated by contractions of edges
(that is, there is a morphism from $\tau$ to $\tau'$ if $\tau'$ is the contraction
of $\tau$ by $e$) and relabelings which are isomorphisms
(that is, there is a morphism from $\tau$
to $\tau'$, which could be $\tau$ itself, if $\tau'$ is obtained from $\tau$ by
relabeling of its edges).
\end{itemize}
\end{definition}

See \refF{basicmor} for some examples of objects and morphisms in $\Upsilon$.
Let  $\Upsilon_n$ denote the full subcategory
of trees with $n$ leaves, which has a terminal object, namely the unique tree with one vertex
called the $n$th corolla $\gamma_n$.  We allow for the tree
$\gamma_0$ which has no leaves, only a root vertex, and is
the only element of $\Upsilon_0$.
For a vertex $v$ let $|v|$ denote the number of edges for which
$v$ is terminal, usually called the arity of $v$.

\begin{definition}\label{D:op}
An operad is a functor $\Op$ from $\Upsilon$ to a symmetric
monoidal category $(\mC, \odot)$ which satisfies the following axioms.

\begin{enumerate}
\item $\Op(\tau) \cong \odot_{v \in \tau} \Op(\gamma_{|v|})$.  \label{1}
\item If $e$ is a redundant edge and $v$ is its terminal vertex
then $\Op(c_{\{e\}})$ is the identity map on
$\odot_{v' \neq v} F(\gamma_{v'})$ tensored  with the isomorphism
$(\one_{\mC} \odot -)$ under the  decomposition of axiom (\ref{1}).  \label{3}
\item If $S$ is a subtree of $\tau$ and if $f_{S, S'}$ and $f_{\tau, \tau'}$
contract the same set of edges, then under the
decomposition of (\ref{1}), $F(f_{\tau,\tau'}) = F(f_{S, S'}) \odot id$. \label{4}
\end{enumerate}
\end{definition}

By axiom~(\ref{1}), the values of $\Op$ are determined by 
its values on the corollas $\Op(\gamma_n)$, which corresponds
to $\Op(n)$ in the usual terminology.   Because of the relabeling morphisms,
$\Op(n)$ has an action by  the $n$th symmetric group.  By axiom~(\ref{4}),
the values of $\Op$ on morphisms may be computed by
composing morphisms on sub-trees, so we may identify some 
subset of basic morphisms through which all morphisms
factor.  In \refF{basicmor} we illustrate some basic morphisms in $\Upsilon$.
The first  corresponds to what are known as $\circ_i$ operations.
The second corresponds May's
operad structure maps from Definition~1.1 of \cite{May72}.  
That $\Op$ is a functor implies 
the commutativity of diagrams involving these basic morphisms.

Filling in what the operad structure maps are for our
Examples from ~\ref{Ex:operads} is a pleasant exercise, which we leave
in part to the reader.  
\begin{examples}
\begin{enumerate}
\item For $\com$, all structure maps are the identity.
\item For $\ass$, they are ``insertion and relabeling.''
\item For $\lie$, the structure maps are defined
by grafting trees, that is identifying the root of one with the leaf
of another.  These are well-defined because the Jacobi and anti-symmetry
identities are defined locally.  
\item  For $\dis^{d}$ we give a full account.
Let $T$ be a tree whose vertices consist of the root
vertex $v_0$ and a terminal vertex $v_e$ for each root edge $e$.  
Thus, $T \to \gamma_n$, where $n$ is the number of leaves of $T$,
gives rise to one of May's structure maps as in \refF{basicmor}.  
Given a label $i \in \n$ let $v(i)$ be the initial vertex for the $i$th
leaf, let $o(i)$ be the label of leaf $i$ within the ordering on edges
of $v(i)$ and let $e(i)$ be the label of the root edge for
which $v(i)$ is terminal.
Define $\dis^{d}(T \to \gamma_n)$ as follows 
$$(x^v_i, r^v_i)^{v \in V(T)}_{1 \leq i \leq \#v} \mapsto (y_j, \rho_j)_{j \in \n} \;\;\; {\rm where}
\;\;\; y_j = x^{v_0}_{e(j)} + r^{v_0}_{e(j)}  x^{v(j)}_{o(j)} \;\;\; 
{\rm and} \;\;\; \rho_j = r^{v_0}_{e(j)} r^{v(j)}_{o(j)}.$$
See Figure~6 for the standard picture.
\begin{figure}[ht]\label{F:disks}
$$\includegraphics[width=12cm]{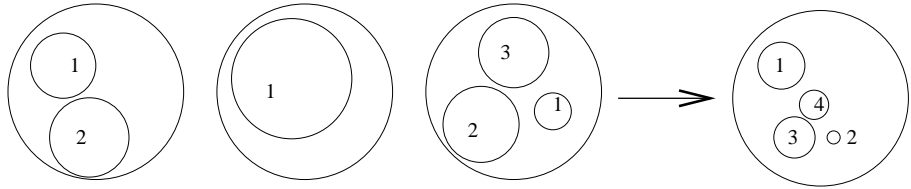}$$ 
\caption{A structure map for the $2$-disks operad.} 
 \end{figure}
\item In the case of $\poid$, the structure maps are essentially grafting as for 
$\lie$, but with an important additional wrinkle given by the Leibniz rule.
In order to be precise without unnecessary complication, it helps to switch
from forests to more algebraic notation.  We may associate to an $n$-forest 
an expression in variables $x_{1}, \cdots, x_{n}$ with two binary products,
denoted $\cdot$ and $[\;,\;]$.  For example to the forest $F =
\begin{xy}
  (1.5,1.5); (3,3)**\dir{-}, 
  (0,3); (3,0)**\dir{-}; 
  (7.5,4.5)**\dir{-},        
  (3,0); (3,-1.5)**\dir{-}, 
  (3,6); (6,3)**\dir{-},
  (6,6); (4.5,4.5)**\dir{-},
  (0,4.2)*{\scriptstyle 2},
  (3,4.2)*{\scriptstyle 6},
  (3,7.2)*{\scriptstyle 1},
  (6,7.2)*{\scriptstyle 7},
  (7.5,5.7)*{\scriptstyle 3},
\end{xy}
 \begin{xy}                          %
   (0,1.5); (1.5,0)**\dir{-};            
   (3,1.5)**\dir{-},                   
   (1.5,-1.5); (1.5,0)**\dir{-},           
   (0,3.3)*{\scriptstyle 4},        %
   (3,3.3)*{\scriptstyle 5},        %
 \end{xy}$ we associate the bracket expression $\left[[x_{2}, x_{6}], [[x_{1}, x_{7}], x_{3}]\right] \cdot 
 [x_{4}, x_{5}]$.  More generally, bracket expressions may include multiplications
 by $\cdot$ within brackets, but these may be reduced to expressions associated to forests
 after the Leibniz rule, $[X, Y\cdot Z] = (-1)^{|X||Y|} Y \cdot [X, Z] + [X, Y] \cdot Z$, is imposed.

Let $f: \tau \to \gamma_{n}$ be a morphism in $\Upsilon$ in which $\tau$
is a tree with one internal vertex over the $i$th root edge.  
Define an operad structure on $\poid = \oplus \poid(n)$ by sending $f$
to the map $$\poid(n) \otimes \poid(m) \to \poid(n+m -1)$$ where $B_{1} \otimes B_{2}$
is sent to the bracket expression defined as follows.  The variables in $B_{2}$ are re-labeled
from $x_{i}$ to $x_{m}$.  The variables $x_{j}$ in $B_{1}$ with $j>i$ are re-labeled
by $x_{j+m-1}$.  Finally  $B_{2}$ is substituted for $x_{i}$
in $B_{1}$.  Note that in order to express this in terms of the $n$-forest basis, 
the Leibniz rule would then need to be applied repeatedly.

\item For $\ndo_{X}$, structure maps are defined by composition.
\end{enumerate}
 \end{examples}

Finally, we give an anti-climactic definition of an algebra over an
operad.

\begin{definition}
An algebra structure for $X$ over an operad $\Op$ is a natural transformation
of operads $\Op \to \ndo_{X}$.  
\end{definition}

By adjointness, the maps $\Op(n) \to
{\rm Hom}_{\mC}(X^{\odot n}, X)$ give rise to multiplication maps
$\Op(n) \odot X^{\odot n} \to X$.  Because of the relabeling 
morphisms in $\Upsilon$, these maps are equivariant with respect to the diagonal
symmetric group action on $\Op(n) \odot X^{\odot n}$.

As for examples, algebras over $\com$ are commutative algebras, and similar eponymous results
hold for $\ass$, $\lie$, and $\poid$.  We will discuss the little disks operad 
in the next section.  As for  $\ndo_{X}$, the only general statement is that $X$ is an algebra 
over it.

\medskip

We leave historical remarks about the theory of operads for other papers in this volume.

\section{The homology of the little disks operad}\label{S:diskshomology}

The little disks operad
and its action on iterated loop spaces can trace their lineage to the proof of 
one of the first theorems
in homotopy theory, namely that $\pi_{2}(X)$ is an abelian group.  Stated in our
language, that proof gives a path of possible multiplications of $f$ and $g$, two
elements of the second loop space $\Omega^{2}(X)$, starting at $f \cdot g$ and
ending at $g \cdot f$.  That path lies within the ``little rectangles'' sub-operad of the space
of all multiplications.  If we start with an arbitrary number of maps in any dimension, we are led to 
the little disks action on a $d$-fold loop space.

\begin{definition}
The action of $\dis^{d}$ on $\Omega^{d}(X) = {\rm Map}\left((D^{d}, S^{d-1}),(X, *) \right) $ 
is defined through maps
$$\dis^{d}(n) \times \left(\Omega^{d}(X)^{\times n}\right) \to \Omega^{d}X,$$
which send  $\{ B_{i} \} \times \{ f_{i} \}$ to the map whose restriction to $B_{i}$ is
$f_{i}$ composed with the canonical linear homeomorphism of $B_{i}$ with $D^{d}$.
At points in $D^{d}$ outside of any $B_{i}$, the resulting map is constant at the basepoint.
\end{definition} 

Thus a $d$-fold loop space is an algebra over the little disks operad.  Boardman-Vogt \cite{BoVo73}
and May \cite{May72} showed that the converse is essentially true.  That is if $X$
has an action of the little $d$-disks and $\pi_{0}(X)$, which is necessarily then a monoid,
is in fact a group then $X$ is homotopy equivalent to a $d$-fold loop space.  
This result is known as a ``recognition principle,'' since it gives a criterion for
recognizing iterated loop spaces.

\begin{figure}[ht]\label{F:disksaction}
$$\includegraphics[width=10cm]{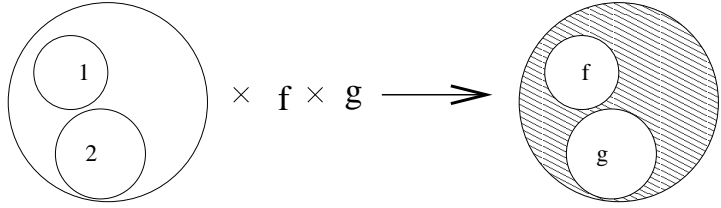}$$ 
\caption{Little $2$-disks acting on two maps.} 
 \end{figure}

In general, operad actions on spaces are reflected in their homology.

\begin{proposition}
The homology of any operad $\Op$ of spaces will
be an operad of modules.  Moreover, the homology of an algebra over
$\Op$ will be an algebra over $H_{*}(\Op)$.
\end{proposition}

\begin{proof}
That $H_{*}(\Op)$ is an operad is immediate by composing the K\"unneth map
$$H_{*}(\Op(r)) \otimes H_{*}(\Op(n_{1})) \otimes \cdots \otimes H_{*}(\Op(n_{r}))
\to H_{*}\left(\Op(r) \times \Op(n_{1}) \times \cdots \times \Op(n_{r})\right)$$
with the induced map in homology of an operad structure map to get a corresponding
operad structure map in homology. Moreover, a map $\Op \to \ndo_{X}$
induces a map $H_{*}(\Op) \to H_{*}(\ndo_{X})$ which in turn maps to 
$\ndo_{H_{*}(X)}$, again using the K\"unneth map.
\end{proof}

We can now state the main result of this paper in full.

\begin{theorem}\label{T:main2}
The homology of the little $d$-disks operad $\dis^{d}$ is 
the degree $d$ Poisson operad $\poid$.  Thus the homology of $\Omega^{d}(X)$
is an algebra over $\poid$.
\end{theorem}

Before establishing this theorem, we reflect on its significance, which is to 
endow the homology of iterated loop spaces a rich additional structure.
This homology has intrinsic interest, but may also be used to study homotopy groups.
The standard Hurewicz map $\pi_{n}(X) \to H_{n}(X)$ is often zero (for example,
in all but one degree when $X$ is  a sphere).  But if we use 
adjointess to identify $\pi_{n}(X)$ with $\pi_{n-k}(\Omega^{k} X)$ then 
the $k$-looped Hurewicz map to $H_{n-k}(\Omega^{k} X)$ can yield
additional information.  For example, a theorem of Milnor and Moore 
(whose proof at the end of \cite{MiMo65} is left as a nice exercise)
states that for rational homotopy and homology, the $1$-looped Hurewicz map is
injective for simply connected spaces.

We now bring in what we know about configuration spaces
through the following.

\begin{lemma}\label{L:equiv}
The space $\dis^{d}(n)$ is homotopy equivalent to $\Co_{n}(\R^{d})$.
\end{lemma}

\begin{proof}
Because ${\rm Int}D^{d}$ is homeomorphic to $\R^{d}$, their associated
configuration spaces $\Co_{n}({\rm Int}D^{d})$ and $\Co_{n}(\R^{d})$
are homeomorphic.  The space of little disks $\dis^{d}(n)$ projects onto
$\Co_{n}({\rm Int}D^{d})$ by definition (mapping to the configuration defined
by the centers of the disks).  This projection defines a fiber bundle whose fibers,
given by the set of possible radii, are convex spaces and thus contractible.
\end{proof}

\refT{main1} now implies \refT{main2} at the level of underlying vector spaces,
so we may focus on operad structure.
To establish the compatibility between the geometric insertion operad structure of the little
disks and the algebraic insertion operad structure of the Poisson operad is 
easier when considering homology classes
of $\dis^{d}$ represented by trees rather than forests.  The key is to choose
appropriately consistent lifts of the submanifolds $P_{F}$ to the spaces $\dis^{d}(n)$
(now the planets really look like planets, being represented by disks rather
than points).  If for example $F$ is a forest and $T$ is a forest with one component,
a $\circ_{i}$ map would send
$\wt{P}_{F} \times \wt{P}_{T}$ precisely to $\wt{P}_{F'}$, where $F'$ is
the grafting of $T$ onto $F$.   So in homology $F \circ_{i} T$ is the grafting
of $T$ onto the $i$th leaf of $F$ accordingly.  

To manage the general case, we focus on cohomology instead of homology.
It is geometrically easier to establish the linear
dual of \refT{main2}, identifying the cohomology
of $\dis^{d}$ with the cooperad structure on  $\iopd$.
Translating to homology is then a matter of pure algebra and combinatorics,
which is carried out in \cite{Sinh06.2}.  A cooperad is a  
functor from $\Upsilon^{\text{op}}$ to a symmetric monoidal
category which satisfies axioms dual to those of \refD{op}.  Note that 
in associating a dual cooperad to an operad, we are 
not changing the symmetric monoidal product.
For example, a standard cooperad in the category of vector spaces 
is defined using the tensor product rather than the direct sum.

\begin{definition} \label{D:mB}
To an o-tree $\tau$ with $n$ leaves and
two distinct integers $j,k \in \n$ let $v$ be the nadir
of the shortest path between leaves labelled $i$ and $j$ and define
$J_v(j), J_v(k)$ to be the labels of the branches of $v$ over
which leaves $j$ and $k$ lie.

The module $\iopd$, forms a cooperad 
which associates to the morphism $\tau \to \gamma_n$ 
the homomorphism $g_\tau$ sending $G \in \Gamma$ to 
$( {\rm{sign}}\ \pi)^{d} \bigotimes_{v_i} G_{v_i}$.
Here $v_i$ ranges over internal
vertices in $\tau$ and $G_{v_i} \in \Gamma^{|v_i|}$, 
is defined by having for each edge in $G$, say from $j$ to $k$,
an edge from $J_v(j)$ to $J_v(k)$ in $G_v$.
The edges of $G_{v_{i}}$ are ordered in accordance with that 
of the edges in $G$ which give rise to them, and $\pi$ is the
permutation relating this order on all of the edges in 
 $\bigotimes_{v_i} G_{v_i}$ to the ordering within $G$.
\end{definition}

Consider for example when $\tau$ is the first tree from Figure~5, with leaf
labeling given by the planar embedding.
The corresponding cooperad structure map would send a single
graph $G$ on five vertices to the tensor product of two graphs, $G_{r}$ with four
vertices and $G_{v}$ with two vertices.   The graph $G_{v}$ would have an
edge between its two vertices if and only if $G$ had $\linep{3}{4}$ as an edge.
Any edge of $G$ of the form $\linep{1}{3}$ or $\linep{1}{4}$ would give
rise to an edge $\linep{1}{3}$ in $G_{1}$.  The edge $\linep{5}{4}$ in $G$
would give rise to $\linep{4}{3}$ in $G_{1}$.

\begin{theorem}[Thm. 6.8 of  \cite{Sinh06.2}]\label{T:coop}
The cooperad structure on $\iopd$ is linearly dual to that of $\poid$
through the configuration pairing.
\end{theorem}

The key to proving this theorem is that the configuration pairing can be
defined directly on bracket expressions.  Looking at \refD{pairing}, we use innermost pairs
of brackets instead of nadirs of paths 
to define the analogue of the map $\beta_{G,T}$.  Remarkably,
the Leibniz rule is respected by this extended definition of this pairing.
Indeed, the configuration pairing can be viewed as the central algebraic object in this area,
and the anti-symmetry, Jacobi, Leibniz, and Arnold identities arise 
naturally in describing its kernel.

While the operad structure of $\poid$ is more familiar, the cooperad structure
on $\iopd$ is simpler. The operad maps on $\poidn$ require the Leibniz rule
to be applied recursively to reduce to any standard basis, 
while the cooperad maps on $\iopd$ require no such reduction for many standard
bases.  While here we use this cooperad as a useful way to prove a theorem
about the corresponding operad, in work on Koszul duality cooperads play an equal role.

\begin{proof}[Proof of \refT{main2}]
By \refL{equiv} and \refT{main1}, $H_{*}(\dis^d)$ and $\poid$ are isomorphic as vector spaces,
so it suffices to consider their operad structure.  By \refT{perfect} and \refT{coop}, we may instead
establish that the cooperad structures on $\iopd$ and the cohomology
of $\dis^d$ agree.

Let $f: \tau \to \gamma_{n}$ be a morphism in $\Upsilon$, where $\gamma_{n}$ is a corolla,
and let $\prod_{w \in \tau} \dis^{d}(|w|) \to \dis^{d}(n)$ be the corresponding operad
structure map.  Using the ring structure on cohomology, it suffices to understand the pullback
of a generator $a_{ij}$.  By definition, we consider the composite 
 $$\pi_{ij} \circ \dis^{d}(f) : \prod_{w \in \tau} \dis^{d}(|w|) \to \dis^{d}(n) {\rightarrow} S^{d-1}.$$ 
 We apply a homotopy  in which at time $t$ 
 the disks in $\dis^{d}(|w|)$ are all scaled by 
by $t^{h}$ where $h$ is the height of the vertex $w$.
As $t$ approaches zero, $\pi_{ij} \circ \dis^{d}(f)$ approaches projection onto 
the factor of  $\dis^{d}(|v|)$ where $v$ is the nadir
of the shortest path between leaves labelled $i$ and $j$, followed by 
$\pi_{J_{v}(i) J_{v}(j)}$, as in \refD{mB}.  Thus $a_{ij}$ pulls back to $a_{J_{v}(i) J_{v}(j)}$
in the $v$th factor of $\dis^{d}$, in agreement with the cooperad structure on $\iopd$.
\end{proof}

To recap, we have now shown that the homology of a $d$-fold loop space is a Poisson
algebra.  The multiplication is the standard one given by loop-sum.  The bracket, known
as the Browder bracket, reflects
some ``higher commutativity'' of the loop-sum.  If  two homology
classes are represented by (pseudo-)manifolds $M, N \to \Omega^{d} X$, then their
bracket will be represented by the map $S^{d-1} \times M \times N \to \Omega^{d} X$
which when restricted to $v \times M \times N$ ``multiplies $M$ and $N$ in the direction
of $v$.''  
 
 \subsection{Historical notes}
The operad structure on spaces of little disks was determined by Cohen in his thesis 
\cite{CLM76}.  There both the non-equivariant and the much more delicate equivariant
homology of these spaces are determined, though the language of operads is not employed. 
Equivariant homology classes yield operations in the homology of iterated loop spaces \cite{DyLa62},
which are algebras over this operad, which is a main focus of \cite{CLM76} and \cite{May71}.  
The simplest example
is with coefficients modulo two, where for any homology class $x$, we have $[x, x] = 0$.  But
this class can be ``divided by two, using only a hemisphere's worth of 
Browder multiplication.''  The 
result as an operation which ``sends $x$ to $\frac{[x,x]}{2}$.''

\bibliographystyle{amsplain}
\bibliography{references}
\end{document}